\documentclass{amsart}
\usepackage{amsfonts}

\newtheorem{theorem}{Theorem}[section]

\theoremstyle{definition}

\theoremstyle{remark}

\numberwithin{equation}{section}
\theoremstyle{plain}

\copyrightinfo{2001}{enter name of copyright holder}

\begin{document}
\title[Asymptotic tightness ]{A remark on the asymptotic tightness in $%
\ell^{\infty}([a,b])$}
\author{Gane Samb LO}
\dedicatory{Dedicated to the memory of Souleymane Niang, the Senegal
Mathematics School founder}

\begin{abstract}
In this note, we extend a simple criteria for uniform tightness in $%
C(0,1) $, the class of real continuous functions defined on $(0,1)$, given
in Theorem 8.3 of Billingsley to the asymptotic tightness in $%
\ell^{+\infty}([a,b])$, the class of real bounded functions defined on $[a,b]
$ with $a<b$, in the lines of Theorems 1.5.6 and 1.5.7 in van der vaart and
Wellner.
\end{abstract}

\maketitle

\Large

\section{Introduction}

In this note, we adapt a powerful tool of Billingsley \cite{bil}. In order
to describe that tool, we are going to make a number of definitions and
reminders.

\subsection{Uniform tightness}

Let $X_{1},X_{2},..$ be sequence of random elements with values in \ $%
S_1=C(0,1).$ This sequence is said to be tight, that is the sequence of
probability measures $(\mathbb{P}_{X_{n}})_{n\geq 0}$ is tight, if and only if for any $\varepsilon >0$,
there exists a compact set $K_{\varepsilon }$ of $S_1$ such that
\begin{equation*}
\sup_{n\geq 1} \mathbb{P}(X_{n}\in K_{\varepsilon })\leq \varepsilon .
\end{equation*}

\noindent He proved :

\begin{theorem}
\label{b1} \bigskip The sequence of the probability measures $(\mathbb{P}_{X_{n}})_{n\geq 0}$ is tight in $S_{1}$ if and only if

\begin{itemize}
\item[(i)] The sequence of the probability measures $(\mathbb{P}
_{X_{n}(0)})_{n\geq 0}$ is tight in $\mathbb{R}$ and

\item[(ii)] The sequence $X_{n}$ is uniformly equicontinuous in
probability, that is, for any $\eta >0,$%
\begin{equation}
\lim_{\delta \rightarrow 0}\sup_{n\geq 1} \mathbb{P}(\sup_{\left\vert
s-t\right\vert <\delta ,(s,t)\in (0,1)^{2}}\left\vert
X_{n}(s)-X_{n}(t)\right\vert >\eta )=0.  \label{ucb}
\end{equation}
\end{itemize}
\end{theorem}

\newpage
\noindent Since (\ref{ucb}) is not easy to handle in general, a stronger criteria is generally used. It is only a sufficient condition for tightness. It is exposed in Billingsley \cite{bil}, as follows

\begin{theorem}
\label{b2} Suppose the two following assertions hold

\begin{itemize}
\item[(i)] The sequence of the probabiliy measures $(\mathbb{P}%
_{X_{n}(0)})_{n\geq 0}$ is tight in $\mathbb{R}$.

\item[(iii)] For or any $\eta >0$,
\begin{equation*}
\lim_{\delta \rightarrow 0} \sup_{s\in (0,1)} \sup_{n\geq 1}\frac{1}{\delta }%
\mathbb{P}(\sup_{s-\delta <t<s+\delta ,t\in (0,1)}\left\vert
X_{n}(s)-X_{n}(t)\right\vert >\eta )=0.
\end{equation*}
\end{itemize}

\bigskip \noindent Then sequence of the probability measures $(\mathbb{P}%
_{X_{n}})_{n\geq 0}$ is tight.
\end{theorem}

\subsection{Asymptotic tightness for non measurable random applications}

Now consider the more general space $S_{2}=\ell ^{\infty }([a,b]),a<b, $ the
set of all bounded and real functions defined on $[a,b]$, equipped with the
supremum norm $\left\Vert x\right\Vert =\sup_{t\in \lbrack a,b]}\left\vert
x(t)\right\vert $. Let $\left( X_{\alpha }\right) _{\alpha\in D}$ be a field
of applications with values in S$_{2}$ and such that each $X_{\alpha}$ is
defined on a probability space $(\Omega_{\alpha}, \mathcal{A}_{\alpha}, 
\mathbb{P}_{\alpha})$ and is not necessarily measurable, where $D$ is a
well-directed set. It is said that the field $\left(X_{\alpha }\right)
_{\alpha \in D}$ is asymptotically measurable if and only if for any real, and bounded and continuous $f$ defined on $S_2$ (denoted $f\in $C$_{b}(S)$), we have
\begin{equation*}
\lim_{\alpha }E^{\ast }f(X_{\alpha })-E_{\ast }f(X_{\alpha })=0.
\end{equation*}

\noindent where $\mathbb{E}^{\ast}$, $\mathbb{E}_{\ast}$, $\mathbb{P}^{\ast}$ and  $\mathbb{P}_{\ast}$ respectively stand for the outer integral, the inner integral, the outer probability and the inner probability.\\

\noindent It is asymptotically tight if and only if for any $\varepsilon >0$, there exists a
compact set $K_{\varepsilon }$ such that for any $\delta >0$

\begin{equation*}
\liminf_{\alpha }\sup \mathbb{P}_{\ast }(X_{n}\in K_{\varepsilon }^{\delta })\geq
1-\varepsilon ,
\end{equation*}

\noindent where $K_{\varepsilon }^{\delta }=\{y\in S_{2},\left\Vert x-K_{\varepsilon
}\right\Vert <\delta \}$ is the $\delta $-enlargement of $K_{\varepsilon}$.\\

\noindent The following characterization of the asymptotic tightness in $\ell
^{\infty }([a,b])$ is given in \cite{var} as follows.

\begin{theorem}
\label{v1} The field $\left( X_{\alpha }\right) _{\alpha \in D}$ is
asymptotically tight if and only if

\begin{itemize}
\item[(iv)] each margin $X_{\alpha }(t)$, $t\in \lbrack a,b \rbrack$, is
asymptotically tight in $\mathbb{R}$ and 

\item[(v)] there exists a semi-metric $\rho $\ on S$_{2}$ such that $%
([a,b],\rho )$ is totally bounded and such that for any \ $\varepsilon $\ $%
>0,$ and  for any \ $\eta >0,$%
\begin{equation}
\lim_{\delta \rightarrow 0}\limsup_{\alpha }{\large \mathbb{P}}^{\ast }%
{\large (}\sup_{\rho (s,t)<\delta }\left\vert X_{\alpha }(s)-X_{\alpha
}(t)\right\vert {\large >\eta )=0.}  \label{ucv}
\end{equation}
\end{itemize}
\end{theorem}

\bigskip \noindent We are going to make comments on these theorems in the next section
where we state the problem and propose a solution.

\section{The Result}

We observe that for $\rho (s,t)=\left\vert s-t\right\vert$, the space $([a,b],\rho )$
is totally bounded and (\ref{ucb}) and (\ref{ucv}) coincide under the
assumption of $a.s.$\ continuity of the $X_{\alpha }.$ Thus, it is natural
to know whether Theorem \ref{v1} has an analogue in $S_{2}.$ Indeed, we have

\begin{theorem}
Assume that the two following assertions hold.

\begin{itemize}
\item[(iv)] Each margin $X_{\alpha }(t)$, $t\in \lbrack a,b],is$
asymptotically tight in $\mathbb{R}$.

\item[(v)] For $s\in \lbrack a,b \rbrack,$ for any $\eta >0,$%
\begin{equation}
\lim_{\delta \rightarrow 0} \sup_{s\in [0,1]}\limsup_{\alpha }\frac{1}{%
\delta }{\large \mathbb{P}}^{\ast }{\large (}\sup_{s-\delta <t<s+\delta
,t\in \lbrack a,b]}\left\vert X_{\alpha }(s)-X_{\alpha }(t)\right\vert 
{\large >\eta )=0.}  \label{ucva}
\end{equation}
\end{itemize}

\bigskip Then the field $\left( X_{\alpha }\right) _{\alpha \in D}$ is
asymptotically tight.
\end{theorem}

\bigskip \noindent \textbf{Proof}. Suppose that $(iv)$ and $(v)$ hold. Let $0<\delta
<(b-a)$ and $\eta >0$. Put 
\begin{equation*}
A_{t}=\{z:\sup_{t\leq s\leq t+\delta }\left\vert z(s)-z(t)\right\vert >\eta).
\end{equation*}

\noindent 
The intervals $I_{i}=[a+i\delta ,a+(i+1)\delta ]$\ form a partition of $%
[a,b]$. Consider a $z\in S_{2}$ such that%
\begin{equation}
\forall i\leq \frac{b-a}{\delta },z\notin A_{t_{i}}  \label{cc}
\end{equation}%
where $t_{i}=a+i\delta \leq b.$ Let $(s,t)\in \lbrack a,b]$ such that $%
\left\vert s-t\right\vert <\delta $. Then either $s$ and $t$ lie in the same
interval $I_{i}$ or they lie in adjacent ones. In the latter case, put $t\in I_{i}$\ et $s\in I_{i+1}$, where we suppose that $t\leq s.$ We have $%
t_{i}=a+i\delta $ and 
\begin{equation*}
(\ref{cc})\Rightarrow \left\vert z(s)-z(t)\right\vert \leq \left\vert
z(s)-z(t_{i})\right\vert +\left\vert z(t_{i})-z(t_{i+1})\right\vert
+\left\vert z(t_{i+1})-z(t)\right\vert <3\eta .
\end{equation*}

\noindent This implies. 
\begin{equation*}
\sup_{\left\vert s-t\right\vert <\delta }\left\vert z(s)-z(t)\right\vert
\leq 2\eta .
\end{equation*}

\noindent We get that 
\begin{equation*}
z\in \{x:\sup_{\left\vert s-t\right\vert <\delta }\left\vert
x(s)-x(t)\right\vert \geq 3\eta )
\end{equation*}

\noindent implies that there exists an indice $i$ such that 
\begin{equation*}
z\in A_{t_{i}}.
\end{equation*}

\noindent Then 
\begin{equation*}
\{z:\sup_{\left\vert s-t\right\vert <\delta }\left\vert z(s)-z(t)\right\vert
\geq 3\eta )\subset \bigcup_{i}A_{t_{i}}.
\end{equation*}

\noindent Hence 
\begin{equation*}
\mathbb{P}^{\ast }(X_{\alpha }\in \{z:\sup_{\left\vert s-t\right\vert
<\delta }\left\vert z(s)-z(t)\right\vert \geq 3\eta \})\leq \sum_{i\leq
(b-a)/\delta }\mathbb{P}^{\ast }(X_{\alpha }\in A_{t_{i}}).
\end{equation*}

\noindent Thus 
\begin{equation*}
\mathbb{P}^{\ast }(\sup_{\left\vert s-t\right\vert <\delta }\left\vert
X_{\alpha }(s)-X_{\alpha }(t)\right\vert \geq 3\eta )\leq \sum_{i\leq
(b-a)/\delta }\mathbb{P}^{\ast }(\sup_{t_{i}\leq s\leq t_{i}+\delta
}\left\vert X_{\alpha }(s)-X_{\alpha }(t)\right\vert >\eta ).
\end{equation*}

\noindent We apply (\ref{ucva}) with $3\eta $ to get 
\begin{equation*}
\limsup_{\alpha }\mathbb{P}^{\ast}(\sup_{\left\vert s-t\right\vert
<\delta }\left\vert X_{\alpha }(s)-X_{\alpha }(t)\right\vert \geq 3\eta )
\end{equation*}

\begin{equation}
\leq \limsup_{\alpha }\leq \delta \sum_{i\leq (b-a)/\delta }\delta ^{-1}\mathbb{P}^{\ast }(\sup_{t_{i}\leq s\leq t_{i}+\delta }\left\vert X_{\alpha
}(s)-X_{\alpha }(t)\right\vert >3\eta )  \notag
\end{equation}%
\begin{equation*}
\leq \delta ([\frac{b-a}{\delta }]+1)\limsup_{\alpha }\max_{i\leq
(b-a)/\delta }\left( \delta ^{-1}\mathbb{P}^{\ast }(\sup_{t_{i}\leq s\leq
t_{i}+\delta }\left\vert X_{\alpha }(s)-X_{\alpha }(t)\right\vert >3\eta
)\right) .
\end{equation*}

\noindent It comes that 
\begin{equation*}
\lim_{\delta \rightarrow 0}\limsup_{\alpha }\mathbb{P}^{\ast
}(\sup_{\left\vert s-t\right\vert <\delta }\left\vert X_{\alpha
}(s)-X_{\alpha }(t)\right\vert \geq 3\eta )=0.
\end{equation*}

\noindent The proof is complete.

\bigskip

\end{document}